\documentclass{amsart}
\usepackage{amssymb} 
\newcommand\datver[1]{\def\datverp 
 {\par\text{Version: #1; Run: \today}}} 
\datver{3.0A; Revised: 6/17/2000} 

\newcommand\CC{\mathbb C} 
 
\newcommand\NN{\mathbb N} 
\newcommand\RR{\mathbb R} 
\newcommand\ZZ{\mathbb Z} 
 
\newcommand\pa{\partial} 
\newcommand\GR{\mathcal{G}} 
\newcommand\CI{\mathcal{C}^\infty} 
\newcommand\CIc{\mathcal{C}_c^\infty} 
\newcommand\Tt{a_0 \otimes a_1 \otimes \ldots \otimes a_n} 
\newcommand\ie{i.e. } 
\newcommand\HH{\operatorname{HH}} 
\newcommand\Hd{\operatorname{HH}} 
\newcommand\Hc{\operatorname{HC}} 
\newcommand\HC{\operatorname{HC}} 
 
\newcommand\Hp{\operatorname{HP}} 
\newcommand\EH{\operatorname{EH}} 
\newcommand\EC{\operatorname{EC}} 
 
\newcommand{\tPS}[1]{\Psi^{#1}(\GR)} 
 
\newcommand\mO{{\mathcal O}(M)} 
\newcommand\mOt{{\mathcal O(M)} \otimes_{\CI(M)}} 
 
\newcommand\bF[2]{{\mathcal M}_{#1}(#2)} 
 
\newcommand{\R}{\RR}

\newcommand{\Z}{\ZZ}

\newcommand{\cohom}{\operatorname{H}} 
\newcommand\alge{\mathcal A} 
\newcommand\ideal{\mathcal I} 
 
\newcommand\maL{\mathcal L} 
\newcommand{\LX}{{\mathcal L}(M)} 
\newcommand\potimes{\widehat{\otimes}} 
 
\newtheorem{theorem}{Theorem} 
\newtheorem{proposition}{Proposition} 
\newtheorem{corollary}{Corollary} 
\newtheorem{lemma}{Lemma} 
 
\theoremstyle{definition} 
\newtheorem{definition}{Definition} 
 
\theoremstyle{remark}

\begin{document} 
 
\title[Homology of complete symbols] 
{Homology of complete symbols and non-commutative 
geometry} 
 
\author[M. Benameur]{Moulay-Tahar Benameur} 
\address{Pennsylvania State University and Inst. Desargues, Lyon, France} 
\email{benameur@desargues.univ-lyon1.fr} 
\author[V. Nistor]{Victor Nistor} 
\address{Inst. Desargues and Pennsylvania State University, 
University Park, PA 16802} 
\email{nistor@math.psu.edu} 
 
\thanks{Partially supported by the NSF Young Investigator Award 
DMS-9457859,  NSF Grant 
DMS-9971951 and "collaborative research." {\bf 
http:{\scriptsize//}www.math.psu.edu{\scriptsize/}nistor{\scriptsize/}}.} 
 
\dedicatory\datverp 

\begin{abstract} 
We identify the periodic cyclic homology of the 
algebra of complete symbols on a differential groupoid $\GR$ in terms
of the cohomology of $S^*(\GR)$, the cosphere bundle of $A(\GR)$,
where $A(\GR)$ is the Lie algebroid of $\GR$.  We also relate the
Hochschild homology of this algebra with the homogeneous Poisson
homology of the space $A^*(\GR) \smallsetminus 0 \cong S^*(\GR) \times
(0,\infty)$, the dual of $A(\GR)$ with the zero section removed. We
use then these results to compute the Hochschild and cyclic homologies
of the algebras of complete symbols associated with manifolds with
corners, when the corresponding Lie algebroid is rationally isomorphic
to the tangent bundle.
\end{abstract} 
 
\maketitle 
\tableofcontents

\section*{Introduction\label{Sec.I}}

Singular cohomology is often used in Algebraic Topology to obtain
invariants of topological spaces. In the same spirit, Hochschild and
cyclic homology often provide interesting invariants of algebras. A
possible important application of these algebra invariants is to the
study of spaces with additional structures; these include, for
instance, spaces with singularities or spaces endowed with group
actions. This is one of the fundamental ideas of non-commutative
geometry; see Connes' book \cite{ConnesBOOK} and the references
therein.
 
Let $\GR$ be a differentiable groupoid with units $M$, a manifold with
corners, and Lie algebroid $A(\GR) \to M$, (see \cite{LauterNistor} in
this volume for definitions, notation, and background material). To
$\GR$ one can associate several algebras:\ the convolution algebras
$\CIc(\GR)$, $L^1(\GR)$, $\tPS{\infty}$, or other variants of these
algebras. These algebras have always been a favorite toy model for
non-commutative geometry and have several applications, see
\cite{BrylinskiNistor,ConnesNCG,ConnesBOOK,ConnesCuntz,Nistorboundary}.
The algebra $\tPS{\infty}$ of pseudodifferential operators on $\GR$,
for example, is expected to play an important role in the analysis on
singular spaces \cite{LauterNistor}. The ideal $\tPS{-\infty}$ of
regularizing operators in $\tPS{\infty}$ identifies with $\CIc(\GR)$
and gives rise to index invariants via the boundary map in algebraic
$K$-theory:
\begin{equation*} 
	\pa : K_1^{alg}(\mathfrak C) \to K_0^{alg}(\tPS{-\infty}), 
\end{equation*} 
where $\mathfrak C = \alge (M):= \tPS{\infty}/\tPS{-\infty}$ or 
$\mathfrak C = \alge_0(M) := \tPS{0}/\tPS{-\infty}$, depending on the 
nature of the problem.  Typically, much more attention has been 
devoted to finding cohomological invariants of $\tPS{-\infty}$, 
because its (algebraic) $K_0$ is so hard to compute. The computation 
of the $K$-theory of algebras associated with groupoids or their 
$C^*$-closures is far from solved and is part of the larger program 
involving the Baum-Connes Conjecture(s) \cite{BaumConnesHigson}. 
 
In this paper, we concentrate on the homology of $\alge (M)$ or 
variants of this algebra.  We are interested in computing the 
Hochschild, the cyclic, and the periodic cyclic homology groups of 
$\alge (M)$, denoted respectively by $\Hd_*(\alge (M))$, by 
$\Hc_*(\alge (M))$, or by $\Hp_*(\alge (M))$. Of all these, the 
periodic cyclic homology is the easiest to compute. The result is in 
terms of $A(\GR)$, the {\em Lie algebroid} of $\GR$.  Let $S^*(\GR)$ 
be the cosphere bundle of $A^*(\GR)$, that is, the set of unit vectors 
in the dual of the Lie algebroid of $\GR$, and denote $\cohom^{[q]} = 
\oplus_{k \in \ZZ} \cohom^{q + 2k}$.

\begin{theorem}\ 
Assume that the base $M$ is $\sigma$-compact, then 
\begin{equation} 
	\Hp_q(\alge (M)) \cong \cohom_c^{[q]}(S^*(\GR) \times S^1)
	\text{ and } \Hp_q(\alge_0 (M) ) \cong
	\cohom_c^{[q]}(S^*(\GR)).
\end{equation} 
\end{theorem}

The assumption that $M$ be $\sigma$-compact can be replaced with the
assumption that $M$ be paracompact, but then we have to work with more
complicated directed sets, and this is usually unnecessary in
practice.
 
Recall that $A^*(\GR)$ has a natural Poisson structure. We do show
that the natural filtration on the complex computing the Hochschild
homology of the algebra $\alge (M)$ gives rise to a spectral sequence
with $E^2$-term identified with the homogeneous Poisson homology
(Definition \ref{poisson}) of $A^*(\GR) \smallsetminus 0$, the dual
vector bundle of $A(\GR)$ with the zero section removed. We expect
this spectral sequence to degenerate at $E^2$ and to be convergent to
the Hochschild homology of $\alge (M)$. The quantization of $A^*(\GR)$
with this Poisson structure was studied in \cite{LandsmanRamazan}.
 
For certain algebras associated to manifolds with corners, we identify
the homogeneous Poisson homology of $A^*(\GR) \smallsetminus 0$ in
terms of a space $\maL(S^*M)$ functorially associated to the base $M$.
Moreover, the particular form of the resulting spectral sequence
guarantees its convergence. This leads to an identification of the
Hochschild homology of the Laurent complete symbols algebra
$\alge_{\maL}(M)$.

\begin{theorem}\ Let $\mathcal{O}(M)$ be the ring of functions with 
only rational (\ie\ Laurent-type) singularities at the faces of $M$. 
Assume that $\mathcal{O}(M) \Gamma(A(\GR)) \cong \mathcal{O}(M) 
\Gamma(TM)$ via the anchor map, then with $n=\dim(M)$ 
\begin{equation} 
	\Hd_q(\alge_{\maL}(M)) \cong \cohom_c^{2n-q}(\maL(S^*M) \times
	S^1).
\end{equation} 
\end{theorem}

A similar result holds true in the relative case of symbols vanishing 
to infinite order at some subset of $M$, thus extending results of 
\cite{MelroseNistor2} from the case of manifolds with boundary to that 
of manifolds with corners. 
 
In \cite{MelroseNistor1}, the norm closure of the algebra of 
pseudodifferential operators on a manifold with corners was studied 
from the point of view of $K$-theory. However, the $K$-theory is 
sometimes too rough to identify more subtle invariants -- like the 
$\eta$-invariant of Atiyah, Patodi, and Singer 
\cite{AtiyahPatodiSinger1} -- that are not homotopy invariant. This 
was partly remedied in \cite{MelroseNistor1}, where the Hochschild 
$1$-cocycle that gives the index was identified in terms of 
residues. This cocycle was then split into two parts that are direct 
analogues of the Atiyah-Singer integrand and, respectively, the 
$\eta$-invariant. 
 
We now describe briefly the contents of each section.  In Section 
\ref{Sec.HandC}, we introduce the class of algebras we shall work with, 
that is, the class of ``topologically filtered algebras,'' (Definition 
\ref{def.t.filtered}), a class of algebras for which the 
multiplication is not jointly continuous, but which still has a weak
continuity property for multiplication. Because of this, the usual
definitions of the Hochschild and cyclic complexes of a topologically
filtered algebra have to be adapted to our more general
framework. Namely, we have to use iterated inductive and projective
limits. Then we establish some results on the spectral sequences
associated to the natural filtrations of the resulting complexes. In
Section \ref{Sec.PREM}, we establish some technical results on de Rham
complexes with singularities for manifolds with corners, in the spirit
of \cite{MelroseNistor2}. In Section \ref{Sec.gen.HOM} we identify the
periodic cyclic cohomology of the algebra of complete symbols and
relate the Hochschild homology of those algebras with the
(homogeneous) Poisson cohomology of $A^*(\GR) \setminus 0$.  In
Section \ref{Sec.HH}, we compute the Hochschild homology of the
algebra of complete symbols when $A(\GR)$ is rationally isomorphic to
$TM$.  For manifolds without corners (or boundary), these results are
due to Wodzicki and Brylinski-Getzler \cite{BrylinskiGetzler}. Other
related results were obtained by Lauter-Moroianu \cite{LauterMoroianu}
and Moroianu \cite{Moroianu}. We then use these results in Section
\ref{Sec.Applications} to study residues and to determine the cyclic
homology of the algebra of complete symbols, still assuming that
$A(\GR)$ is rationally isomorphic to $TM$. The appendix contains a
short review of projective and inductive limits.

\subsubsection*{Acknowledgements}\ We would like to thank 
Alain Connes, Thierry Fack, Robert Lauter, Sergiu Moroianu, and Serge 
Parmentier for useful discussions.

\section{Hochschild and cyclic homology of filtered 
algebras\label{Sec.HandC}}

We begin by recalling the definitions of Hochschild and cyclic 
homology groups of a topological algebra $\alge$. A good reference is 
Connes' book \cite{ConnesBOOK}. These definitions have to be 
(slightly) modified when the multiplication of our algebra is only 
separately continuous. We thus discuss also the changes necessary to handle 
the class of algebras we are interested in, that of ``topologically 
filtered algebras'' (Definition \ref{def.t.filtered}), and then we 
prove some results on the homology of these algebras. 
 
First we consider the case of a topological algebra $\alge$. Here 
``topological algebra'' has the usual meaning, $\alge$ is a real or 
complex algebra, which is at the same time a locally convex space such 
that the multiplication $\alge \times \alge \to \alge$ is 
continuous. Denote by $\potimes$ the projective tensor product and 
$\mathcal H_n(\alge):= \alge^{\potimes n +1}$, the completion of 
$\alge^{\otimes n+1}$ in the topology of the projective tensor 
product. Also, we denote as usual by $b'$ and $b$ the Hochschild 
differentials: 
\begin{equation}\label{eq.def.bb'} 
\begin{gathered} 
	b'(\Tt)=\sum_{i=0}^{n-1} (-1)^ia_0\otimes\ldots\otimes a_i
	a_{i+1}\otimes\ldots\otimes a_n,\\ b(\Tt)=b'(\Tt) +(-1)^n
	a_na_0\otimes\ldots\otimes a_{n-1}.
\end{gathered} 
\end{equation} 
 
The {\em Hochschild homology groups} of the topological algebra 
$\alge$, denoted $\Hd_*(\alge)$, are then the homology groups of the 
complex $({\mathcal H}_n(\alge), b).$ By contrast, the complex 
$({\mathcal H}_n(\alge), b')$ is often acyclic, for example when 
$\alge$ has a unit. A topological algebra $\alge$ for which 
$({\mathcal H}_n(\alge), b')$ is acyclic is called {\em $H$-unital} 
(or, better, {\em topological $H$-unital}), following Wodzicki 
\cite{Wodzicki}. 
 
We now define cyclic homology. We shall use the notation of 
\cite{ConnesNCG}: 
\begin{equation} 
\begin{gathered} 
	s(\Tt)=1\otimes \Tt, \\ t(\Tt)=(-1)^n a_n\otimes
	a_0\otimes\ldots\otimes a_{n-1}, \\ B_0(\Tt)=s\sum_{k=0}^{n}
	t^k(\Tt), \;\; \text{ and } \; B = (1 - t)B_0.
\end{gathered} 
\end{equation} 
Then $[b,B]_+ := bB + Bb = B^2 = b^2 = 0$, and hence, if we define 
\begin{equation} 
	{\mathcal C}(\alge)_n=\bigoplus_{k\geq 0} {\mathcal H}_{n -
	2k}(\alge),
\end{equation} 
$({\mathcal C}(\alge),b+B),$ is a complex, called {\em the cyclic 
complex} of $\alge$, whose homology is by definition the {\em cyclic 
homology} of $\alge$, as introduced in \cite{ConnesNCG} and 
\cite{Tsygan}. 
 
Consideration of the natural periodicity morphism $\mathcal C_n(\alge) 
\to \mathcal C_{n-2}(\alge)$ easily shows that cyclic and Hochschild 
homology are related by a long exact sequence 
\begin{equation} \label{eq.SBI} 
	\ldots\rightarrow \Hd_{n}(\alge)\stackrel{I}{\longrightarrow}
	\Hc_{n}(\alge)\stackrel{S}{\longrightarrow}
	\Hc_{n-2}(\alge)\stackrel{B}{\longrightarrow}
	\Hd_{n-1}(\alge)\stackrel{I}{\longrightarrow}\ldots\, ,
\end{equation} 
with the maps $I$, $B$, and $S$ explicitly determined. The map $S$ is 
also called the {\em periodicity} operator.  See 
\cite{ConnesNCG,Loday-Quillen1}. This exact sequence exists whether or 
not $\alge$ is endowed with a topology. 
 
For the algebras that we are interested in, however, the 
multiplication is usually only separately continuous, but there will 
exist an increasing multi-filtration $F_{n,l}^m \subset \alge$ of 
$\alge$, 
\begin{equation*} 
	F_{n,l}^m \alge \subset F_{n',l'}^{m'}\alge, \quad \text{if }
	n \le n',\, l\le l',\, \text{ and } m \le m',
\end{equation*} 
by closed subspaces satisfying the following properties: 
\begin{enumerate} 
\item\ $\alge = \displaystyle{ \cup_{n,l,m} } F^m_{n,l} \alge$; 
\item\ The union $\alge_n := \displaystyle{ \cup_{m,l} F_{n,l}^m \alge }$ is 
a closed 
subspace such that 
\begin{equation*} 
	F_{n,l}^m \alge = \alge_n \cap \big ( \cup_{n} F_{n,l}^m \alge 
	\big ); 
\end{equation*}
\item\ Multiplication maps $F_{n,l}^m \alge \otimes F_{n',l'}^m \alge$ 
to $F_{n + n',l+l'}^m \alge$; 
\item\ The maps $ F_{n,l}^m\alge /F_{n-1,l}^m \alge \otimes 
F_{n',l'}^m \alge /F_{n'-1,l'}^m \alge \to F_{n+n',l+l'}^m\alge /F_{n 
+ n' -1, l+l'}^m \alge$ induced by multiplication are continuous; 
\item\ The quotient $F_{n,l}^m\alge/ F_{n-1,l}^m\alge$ is a nuclear  Frechet 
space in the induced topology; 
\item\ The natural map 
$$ 
F_{n,l}^m \alge \to \displaystyle{\lim_{\leftarrow}}\, 
F_{n,l}^m / F_{n-j,l}^m \alge, \quad j \to \infty 
$$ 
is a homeomorphism; and 
\item\ The topology on $\alge$ is the strict inductive limit of the 
subspaces $F_{n,n}^n\alge$, as $n \to \infty$ (recall that 
$F_{n,n}^n\alge$ is closed in $F_{n+1,n+1}^{n+1}\alge$). 
\end{enumerate}

\begin{definition}\label{def.t.filtered}\ 
An algebra $\alge$ satisfying the conditions 1--7, above, will be
called a {\em topologically filtered} algebra.
\end{definition}

It follows from the definition that if $\alge^m :=\displaystyle{ 
\cup_{n,l} F^m_{n,l} \alge }$, then $\alge^m$ is actually a subalgebra 
of $\alge$ which is topologically filtered in its own, but with 
multi-filtration independent of $m$.

For topologically filtered algebras, the multiplication is not 
necessarily continuous, and the definition of the Hochschild and 
cyclic homologies using the projective tensor product of the algebra 
$\alge$ with itself, as above, does not make much sense. For this 
reason, we change the definition of the space $\mathcal H_m(\alge)$ to 
be an inductive limit: 
\begin{equation*} 
	\mathcal H_q(\alge) \cong \displaystyle{\lim_\to} (F_{n,n}^n
	\alge)^{\potimes q+1},
\end{equation*} 
the tensor product being the (complete) projective tensor product. 
The Hochschild homology of $\alge$ is then still the homology of the 
complex $(\mathcal H(\alge), b)$. Since the projective tensor product 
is compatible with the projective limits, we also have 
\begin{equation} 
	\mathcal H_q(\alge) = \lim_{\to} \left( \lim_{\leftarrow}\,
	(F_{n,n}^n \alge / F_{k,n}^n \alge)^{\potimes q+1} \right ),
\end{equation} 
with the induced topology, where first $k \to -\infty$ (in the 
projective limit) and then $n \to \infty$ (in the inductive limit). 
The operator $B$ extends to a well defined map $B : \mathcal 
H_q(\alge) \to \mathcal H_{q + 1}(\alge)$, which allows us to define 
the cyclic complex and the cyclic homology of the algebra $\alge$ as 
the homology of the complex $(\mathcal C_*(\alge),b+B)$, with 
$\mathcal C_q(\alge):= \oplus \mathcal H_{q - 2k}(\alge)$, as for 
topological algebras. 
 
We also observe that both the Hochschild and cyclic complexes have 
natural filtrations given by 
\begin{equation}\label{eq.def.Hfiltr} 
	F_p\mathcal H_q(\alge) := \displaystyle{\lim_{\to}} \left(
	\displaystyle{\lim_{\leftarrow}} \widehat{\otimes}_{j=0}^q \,
	\big( F_{k_j,m}^m \alge / F_{k,m}^m \alge \big ) \right ),
\end{equation} 
where $k_0 + \ldots + k_n \leq p$ defines the filtration. The 
projective and inductive limits are such that first $k,l \to - \infty$ 
and then $m \to \infty$. 
 
For any topologically filtered algebra, we denote 
\begin{equation*} 
	Gr(\alge) := \oplus_n \, \alge_n/\alge_{n-1} 
\end{equation*} 
the {\em graded algebra} associated to $\alge$ (recall that $\alge_n$ 
was defined as $\cup_{l,m}F_{n,l}^m \alge$). Its topology is that of 
an inductive limit of Frechet spaces: 
\begin{equation*} 
	Gr(\alge) \cong \displaystyle{\lim_{N,m,l \to \infty}}
	\oplus_{n = -N}^{N} F_{n,l}^m\alge/F_{n-1,l}^m\alge.
\end{equation*} 
For the algebras like $Gr(\alge)$, we need yet a third way of 
topologizing its iterated tensor products. The correct definition is 
\begin{equation*} 
	\mathcal H_q(Gr(\alge)) \cong \displaystyle{\lim_{N, m ,l \to
	\infty}} (\oplus_{n = -N}^{N}
	F_{n,l}^m\alge/F_{n-1,l}^m\alge)^{\potimes q+1}.
\end{equation*} 
The Hochschild homology of $Gr(\alge)$ is the homology of the complex 
$(\mathcal H_*(Gr(\alge))$,$b)$. The operator $B$ again extends to a 
map $B : \mathcal H_q(Gr(\alge)) \to \mathcal H_{q + 
1}(Gr(\alge))$ and we can define the cyclic homology of $Gr(\alge)$ as 
above. 
 
The Hochschild and cyclic complexes of the algebra $Gr(\alge)$ 
decompose naturally as direct sums of complexes indexed by $p \in 
\ZZ$. For example, $\mathcal H_q(Gr(\alge))$ is the direct sum of the 
subspaces $\mathcal H_q(Gr(\alge))_p$, where 
\begin{multline*} 
	\mathcal{H}_{q}(Gr(\alge))_p = \lim_{m, N ,l \to \infty}
	\bigoplus_{k_j} \, \left ( \widehat{\otimes}_{j=0}^n
	F^m_{k_j,l}\alge/F^m_{k_j-1,l}\alge \right ), \\ \text{ where
	}\; k_0 + k_1 + \ldots + k_n = p \, \text{ and } -N \le k_j
	\le N,
\end{multline*} 
with the induced topology. The corresponding subcomplexes of the 
cyclic complex are defined similarly. We denote by 
$\HH_*(Gr(\alge))_p$ and $\HC_*(Gr(\alge))_p$ the homologies of the 
corresponding complexes.

\begin{lemma}\label{lemma.sp.sq}\ Let $\alge$ be a topologically 
filtered algebra.  Then the natural filtrations on the Hochschild and
cyclic complexes of $\alge$ define spectral sequences $\EH^r_{k,h}$
and $\EC^r_{k,h}$ such that
\begin{equation*}
	 \EH^1_{k,h} \simeq \HH_{k+h}(Gr(\alge))_{k} \;\; \text{and }
	 \EC^1_{k,h} \simeq \HC_{k+h}(Gr(\alge))_{k}.
\end{equation*}
Moreover, the periodicity morphism $S$ induces a morphism of spectral
sequences $S : \EC^r_{k,h} \to
\EC^r_{k,h-2}$, which for $r = 1$ is the graded map associated to the
periodicity operator $S:\HC_{n}(Gr(\alge)) \to \HC_{n-2}(Gr(\alge))$.
\end{lemma}

\begin{proof}\ 
We shall write $F_p = F_p\mathcal{H}(\alge)$, for simplicity, where
the filtration is as defined in Equation \eqref{eq.def.Hfiltr}.
 
The filtration of the complex computing the Hochschild homology of
$\alge$ then gives rise to a spectral sequence with $E^1_{k,h} =
H_{k+h}(F_k/F_{k-1})$, by standard homological algebra. By the
definition of the Hochschild complex of $Gr(\alge)$, 
\begin{equation*}
	\cohom_{k+h}(F_k/F_{k-1}) \cong \HH_{k+h}(Gr(\alge))_{k}. 
\end{equation*} 
This completes the proof for Hochschild homology. The proof for cyclic 
homology is similar 
\end{proof}

\begin{lemma}\label{lemma.conv1}\ 
Let $N$ and $M$ be integers. If $\alge$ is a topologically filtered
algebra with $F_{n,l}^m \alge$ independent of $m$ and $\EH^{1}_{k,h} =
0$ for all pairs $(k,h)$ such that $k < N$ and $k + h \ge M$, then
\begin{equation*}
	\HH_{q}(\alge) \cong \oplus_{k = N}^\infty \EH^{\infty}_{k,q -
	k},
\end{equation*}
if $q \ge M$.  A similar result is true for the cyclic homology 
spectral sequence. 
\end{lemma}

The above isomorphism is not natural, in general, but comes from a
filtration of $\HH_q(\alge)$ whose subquotients identify naturally
with $\EH^{\infty}_{k,q-k}$, see \cite{MacLane,MacLaneMoerdijk}.  What
the above lemma says, put differently, is that the spectral sequence
$\EH^r_{k,h}$ converges to $\HH_{k+h}(\alge)$, for $k + h \ge M$ .

\begin{proof}\ 
We shall use 
\begin{equation*}
	\mathcal{H}_q(\alge) = \displaystyle{\lim\limits_{\to}} \big
	(\lim_{\leftarrow}\, F_{p}\mathcal{H}_q(\alge)/ F_{p'}
	\mathcal{H}_q(\alge) \big).
\end{equation*} 
Denote $F_r=F_{r}\mathcal{H}(\alge)$ for simplicity.  The projective 
limits give rise to a $\lim^1$ exact sequence 
\begin{equation*} 
	0\longrightarrow {\lim_{\leftarrow}}^1
	\cohom_{m+1}(F_p/F_{p'}) \longrightarrow {\cohom}_{m}(F_p)
	\longrightarrow \lim_{\leftarrow}\, {\cohom}_{m}(F_p/F_{p'})
	\longrightarrow 0
\end{equation*} 
for every fixed $p$ (see Lemma \ref{A.exact.limp} from the 
Appendix). {}From the assumption that ${\cohom}_{q}(F_{p'}/F_{p'-1}) 
=0$ for $p' < N$ and $q \ge M$, we know that 
${\cohom}_{q}(F_p/F_{p'})$ becomes stationary for $p' < N$ and $q \ge 
M$. This shows that the $\lim^1$ term above vanishes for $q \ge M$, 
and hence 
\begin{equation*}
	{\cohom}_{q}(F_p) \simeq \displaystyle{\lim_{\leftarrow}}\,
	{\cohom}_{q}(F_p/F_{p'})= {\cohom}_{q}(F_p/F_{N-1}),
\end{equation*}
if $q \ge M$. It also gives that the natural morphism 
$\mathcal{H}(\alge) \to \mathcal{H}(\alge)/F_{N-1}$ induces an 
isomorphism of the $E^r_{k,h}$-terms of the corresponding spectral 
sequences, for $r \ge 1$ and $k + h \ge M$. 
 
Using then the fact that homology and {\em inductive} limits commute 
we obtain 
\begin{multline*} 
	\Hd_q(\alge) = {\cohom}_q (\lim_{p \to
	\infty}\lim_{\leftarrow}\, F_p/F_{p'}) \cong
	\displaystyle{\lim_{p \to \infty}} {\cohom}_q
	(\lim_{\leftarrow} F_p/F_{p'}) \\ \cong \displaystyle{\lim_{p
	\to \infty}} {\cohom}_q (F_p/F_{N-1}) \cong {\cohom}_q
	(\mathcal{H}(\alge)/F_{N-1}) \cong \oplus_{l=N}^\infty
	E^\infty_{k,q - k},
\end{multline*} 
$p' \to -\infty$, where for the last isomorphism we have used that the 
spectral sequence associated to $\mathcal{H}(\alge)/F_{N-1}$ is 
convergent because $\mathcal{H}(\alge)/F_{N-1} = \cup_q F_q/F_{N-1}$. 
\end{proof}

\begin{theorem}\label{theorem.conv1}\ 
Fix an integer $N$ and $M$. If $\alge$ is a topologically filtered
algebra such that each of the algebras $\alge^m := \cup_{n,l}
F_{n,l}^m\alge$ satisfies the assumptions of Lemma \ref{lemma.conv1}
for the given $N$ and $M$, then 
\begin{equation*} 
	\HH_{q}(\alge) \cong \oplus_{k = N}^\infty \EH^{\infty}_{k,q -
	k}, \quad q \ge M.
\end{equation*}
A similar result is true for the cyclic homology spectral sequence. 
\end{theorem}

\begin{proof} \
Each of the algebras $\alge^m = \cup_{n,l}F_{n,l}^m \alge $ is a topologically
filtered algebra in its own if we let $F_{n',l'}^{m'} \alge^m=
F_{n',l'}^m \alge$ (so the filtration of $\alge^m$ 
depends only on $n$ and $l$).
 
For the algebras $\alge^m$, the Hochschild complex is defined as for 
any topologically filtered algebra, except that there is no need to 
take an additional direct limit with respect to $m$. 
 
The Hochschild complex of $\alge$ is then given by 
\begin{equation*}
	\mathcal H_q(\alge) = \lim_{\to} \mathcal H_q(\alge^m), \quad
	m \to \infty
\end{equation*} 
Because taking homology is compatible with inductive limits, we obtain 
that 
\begin{equation*} 
	\HH_q(\alge) = \lim_{\to} \HH_q(\alge^m), \quad m \to \infty. 
\end{equation*} 
Denote by $\EH^r_{k,h}(\alge^m)$ the spectral sequence associated by
Lemma \ref{lemma.sp.sq} to the topologically filtered algebra
$\alge^m$. Again because homology is compatible with inductive limits,
$\EH^r_{k,h}(\alge) \cong \displaystyle{\lim_{\to}}
\EH^r_{k,h}(\alge^m)$.  The result then follows from 
Lemma~\ref{lemma.conv1}. 
\end{proof}

\begin{lemma}\label{lemma.conv2}\ 
Fix an integer $N$ and $a \ge 1$. If $\alge$ is a topologically 
filtered algebra with $F_{n,l}^m \alge$ independent of $m$ and 
$\EH^{a}_{k,h} = 0$, for all $k < N$, then the spectral sequence 
$\EH^r_{k,h}$ converges to $\HH_{k+h}(\alge)$. More precisely, 
\begin{equation*} 
	\HH_{q}(\alge) \cong \oplus_{k = N}^\infty \EH^{\infty}_{k,q -
	k}.
\end{equation*} 
A similar result is true for the cyclic homology spectral sequence. 
\end{lemma}

\begin{proof}\ 
The assumption that $F_{n,l}^{m} \alge$ is independent of $m$ has as a 
consequence that we need not take inductive limits with respect to $m$ 
in the definition of the Hochschild complex of $\alge$. 
 
Denote $F_p=F_{p}\mathcal{H}(\alge)$ for simplicity.  The Hochschild 
complex of $\alge$ is complete, in the sense that 
\begin{equation*} 
	\mathcal{H}(\alge) = \displaystyle{\lim\limits_{\to}} \big
	(\displaystyle{\lim_{\leftarrow}}\, F_{p}\mathcal{H}(\alge)/
	F_{p'}\mathcal{H}(\alge) \big), \quad p' \to -\infty\, \text{
	and }\, p \to \infty.
\end{equation*} 
The projective limits give rise to a $\lim^1$ exact sequence 
\begin{equation*} 
	0\to {\lim_{\leftarrow}}^1 \cohom_{q+1}(\mathcal 
	H(\alge)/F_{p}) \to \cohom_{q}(\mathcal H(\alge)) \to 
	\lim_{\leftarrow}\, \cohom_{q}(\mathcal H(\alge)/F_{p}) \to 0, 
\end{equation*} 
as $p \to -\infty$, for every fixed $q$ (see Lemma 
\ref{A.exlim1}). 
 
The spectral sequence $E^r_{k,h}(p)$ associated to the complex 
$\mathcal H(\alge)/F_p$ is convergent because $\mathcal H(\alge)/F_p = 
\cup_j F_j/F_p$. Consequently, the homology groups of $\mathcal 
H(\alge)/F_p$ are endowed with a filtration $F_{t}\cohom_q(\mathcal 
H(\alge)/F_p)$ such that 
\begin{equation}\label{eq.filtration.spsq} 
F_{t}\cohom_q(\mathcal H(\alge)/F_p)/ F_{t-1}\cohom_q(\mathcal 
H(\alge)/F_p) \cong E^\infty_{t,q-t}(p) 
\end{equation} 
 
Assume now that $p < N$, the assumption that $E^a_{k,h}=0$ for all $k 
< N$ gives that $E^r_{k,h}(p) = E^r_{k,h}$ for all $k < p < N$ or $p + 
a < k$.  Let $A_n = \cohom_q(\mathcal H(\alge)/F_{N - na})$, $B_n = 
F_{N - na + a} \cohom_q(\mathcal H(\alge)/F_{N - na})$, and $C_n = 
A_n/B_n$, $n \ge 2$. Then equation \eqref{eq.filtration.spsq} gives 
that the natural map $A_{n+1} \to A_n$ gives an isomorphism $C_{n+1} 
\cong C_{n}$ and induces the zero map $B_{n+1} \to B_n$.  Using Lemma 
\ref{A.exact.limp} from the Appendix, we obtain that 
$\displaystyle{\lim_{\leftarrow}}^1 A_n = 0$ and 
$\displaystyle{\lim_{\leftarrow}}\, A_n = C_{n_0}$ for any fixed 
$n_0$. 
 
Because $C_{n_0} = \oplus_{l=N}^\infty E^{\infty}_{l, q-l}$, the 
result follows. 
\end{proof}

\begin{theorem}\label{theorem.conv2}\ 
Fix an integer $N$ and $a \ge 1$. If $\alge$ is a topologically 
filtered algebra such that each of the algebras $\alge^m := \cup_{n,l} 
F_{n,l}^{m}\alge$ satisfies the assumptions of Lemma \ref{lemma.conv2} 
for the given $N$ and $a$, then the spectral sequence $\EH^r_{k,h}$ 
converges to $\HH_{k+h}(\alge)$. More precisely, we have 
\begin{equation*} 
	\HH_{q}(\alge) \cong \oplus_{h = N}^\infty \EH^{\infty}_{k,q -
	k}.
\end{equation*} 
A similar result is true for the cyclic homology spectral sequence. 
\end{theorem}

\begin{proof} 
Each of the algebras $\alge^m$ is a topologically filtered algebra in 
its own, if we let $F_{n',l'}^{m'} \alge^m= F_{n',l'}^{m} \alge$ (so 
the filtration really depends only on $n$ and $l$). We shall write 
$F_{n,l} \alge^m$ instead of $F_{n,l}^{m'} \alge^m$, and hence 
$F_{n,l}^{m} \alge = F_{n,l} \alge^m$.

The Hochschild complex of $\alge$ is then given by 
\begin{equation*} 
\mathcal H_q(\alge) = \lim_{\to} \mathcal H_q(\alge^m), \quad 
m \to \infty 
\end{equation*} 
Because taking homology is compatible with inductive limits, we obtain 
that 
\begin{equation*} 
\HH_q(\alge) = \lim_{\to} \HH_q(\alge^m), \quad m \to \infty. 
\end{equation*} 
Denote by $\EH^r_{k,h}(\alge^m)$ the spectral sequence associated to 
the natural filtration of the Hochschild complex the algebra 
$\alge^m$. Because homology is compatible with inductive limits, 
$\EH^r_{k,h}(\alge) \cong \displaystyle{\lim_{\to}} 
\EH^r_{k,h}(\alge^m)$.  The result then follows from 
Lemma~\ref{lemma.conv2}. 
\end{proof}

We conclude this section with a result that is in the same spirit with 
the above results and was independently observed by Sergiu Moroianu in 
a different setting.

\begin{proposition}\label{prop.Hunital}\ 
If the graded algebra $Gr(\alge)$ of the topologically filtered
algebra $\alge$ is $H$-unital, then $\alge$ is $H$-unital, in the
sense that the complex $(\mathcal H_n(\alge), b')$ is acyclic.
\end{proposition}

\begin{proof}\ 
This is completely analogous to the previous results, so we will be 
sketchy. The natural filtration on the complex $(\mathcal H_n(\alge), 
b')$ induces a spectral sequence whose $E^1$ term is the $b'$-homology 
of $Gr(\alge)$. This spectral sequence is proved to be convergent as 
in any of the above two theorems. 
\end{proof}

\section{A rational Laurent de Rham complex\label{Sec.PREM}}

We now recall some notation and obtain some preliminary results on de 
Rham type complexes used in our computations. Examples of 
topologically filtered algebras will be discussed in the following 
section. 
 
Consider a $\sigma$--compact manifold with corners $M$.  Let $\mO$ be 
the ring of functions on $M$ such that $f \in \mathcal O(M)$ if, and 
only if, on every open subset of $M$ diffeomorphic to $[0,1)^k \times 
\RR^{n-k}$, the function $f$ is of the form $x_1^{-p_1} \ldots 
x_k^{-p_k} g$, with $g$ a smooth function on $M$. (So $\mathcal O(M)$ 
is a quotient ring.) By abuse of terminology, we shall call $\mO$ the 
ring of Laurent rational functions on $M$. Also, we shall call 
$x_1,\ldots,x_k$ the {\em local} defining functions of the hyperfaces 
of $M$ at a point $x$ belonging to a corner of codimension $k$. 
 
If there exists a boundary defining function $\rho_{H}\in\CI(M)$ for 
each hyperface $H$ (that is, $H=\{\rho_H=0\}$, $\rho_H \geq 0$ on $M$, 
and $d\rho_H \not = 0$ on $H$), then 
\begin{equation*} 
	\mO = \rho^{-\ZZ}\CI(M)=\bigcup\limits_{j\in\ZZ}x^j\CI(M), 
\end{equation*} 
where $\rho^{-\ZZ}$ stands for the multiplicative set consisting of
all products $\rho_{H_1}^{n_1}\rho_{H_2}^{n_2}\cdots
\rho_{H_m}^{n_m}$, $n_k \in \ZZ$.  We then say that $M$ has {\em
embedded faces}.  A nice discussion of how to avoid the assumption of
$M$ having embedded faces is due to Monthubert (see
\cite{Monthubert1,Monthubert2} and the references within). It is
essential then to use groupoids.

For each $k$, let $\Omega^k_{{c,\maL}}(M) : = \mOt \Omega_c^k(M)$ be 
the space of $k$-forms with compact support on $M$ with only Laurent 
singularities at the faces.  Recall then that {\em the Laurent--de 
Rham cohomology of} $M$ is the cohomology of the Laurent--de Rham 
complex. 
\begin{equation} 
	\cdots\longrightarrow \Omega_{{c,\maL}}^k(M) \longrightarrow 
	\Omega_{{c,\maL}}^{k+1}(M) \longrightarrow \cdots 
\end{equation} 
with respect to the de Rham differential. We denote by $\cohom_{{ 
\maL}}^*(M)$ these cohomology groups. (These groups were denoted 
$\cohom_b^*(M)$ in \cite{MelroseNistor2}, if $M$ is a manifold with 
boundary.) 
 
We shall use also the complex 
\begin{equation} 
	\cdots\longrightarrow \Omega_{{c,\maL}}^k(M)_x \longrightarrow
	\Omega_{{c,\maL}}^{k+1}(M)_x \longrightarrow \cdots.
\end{equation} 
of germs of the Laurent-de Rham complex at a point $x \in M$.

\begin{lemma}\ 
The cohomology of the complex of stalks of the Laurent-de Rham complex 
at the point $x \in M$ belonging to a corner of codimension $k$ is the 
exterior algebra generated by $d\log x_i$, where $x_i$, $i = 1, 
\ldots, k$, are locally defining functions of the hyperfaces 
containing $x$. 
\end{lemma}

\begin{proof}\ 
This statement is local because we are dealing with germs.  So we can 
assume that $x=(0,0,\ldots,0) \in \RR^n$ and $M=[0,\infty)^k \times 
\RR^{n-k} \subset \RR^n$. The complex whose cohomology we have to 
compute is then the projective tensor product of the corresponding 
complexes for $[0,\infty)$ and $\RR$ an appropriate number of 
times. By the Poincar\'{e} Lemma, the complex corresponding to $\RR$ 
has cohomology only in dimension $0$. The cohomology of germs at $0$ 
of Laurent forms on $[0,\infty)$ is seen to be generated by $1$ in 
dimension $0$ and by $d\log t$ in dimension $1$ ($t >0 $). The 
cohomology of a tensor product of these complexes is isomorphic to the 
tensor product of their cohomologies as graded vector spaces (by the 
topological K\"unneth theorem).  This proves the theorem. 
\end{proof}

We shall denote by $\bF jM$ the set of {\em codimension} $j$ faces of 
$M$.

\begin{theorem}\label{theorem.Lcohom}\ 
Let $M$ be a manifold with corners all of whose hyperfaces $H$ are 
embedded submanifolds with defining function $\rho_H$.  Then the 
Laurent-de Rham cohomology spaces of $M$ can be naturally decomposed 
in terms of the cohomology of its faces as 
\begin{equation} 
	\cohom_{{ \maL}}^k(M) =
	\bigoplus\limits_{j=0}^{k}\bigoplus\limits_{F\in \bF{j}M}
	\cohom^{k-j}(F).
\end{equation} 
\end{theorem}

\begin{proof}\ 
Consider for each face $F$ of $M$ the usual de Rham complex 
\begin{equation} 
	\cdots\longrightarrow \Omega^k(F)
	\stackrel{d}{\longrightarrow} \Omega^{k+1}(F) \longrightarrow
	\cdots
\end{equation} 
whose homology is, by (a variant of) de Rham's theorem, the absolute 
cohomology of $F$. Denote by $H_{1},H_2,\ldots,H_m$ the hyperfaces 
containing $F$ and by $\rho_i$ their defining functions. Fix a local 
product structure in a neighborhood of the face $F$ and choose a 
smooth cutoff function $\phi$ with support in that neighborhood and 
equal to $1$ in a smaller neighborhood of $F$.  The map 
\begin{equation} 
	\Phi_F:\CI(F;\Lambda^k) \ni \alpha \longrightarrow \alpha
	\wedge d(\phi \log{\rho_{1}})\wedge\cdots\wedge
	d(\phi\log{\rho_{m}}) \in \rho^{-\ZZ}\CI(M;\Lambda^{k+m})
\end{equation} 
where the local product decomposition near $F$ is used to lift 
$\alpha$ to a smooth form on $M$, is a chain map. 
 
The de Rham complex is a resolution of the constant sheaf on a 
manifold with corners (no factors of $x_i^{-1}$ are allowed). {}From 
this and the previous lemma we obtain that the cochain map $\oplus 
\Phi_F$, where the sum is take over all faces of $M$, gives an 
isomorphism in cohomology. This proves the proposition. 
\end{proof}

The cohomology of the above complex can be described as the cohomology 
of a space $\LX$ naturally associated to $M$ and defined as 
follows. Consider for each face $F$ of $M$ the space $F \times 
(S^1)^k$, where $k$ is the codimension of the face. Moreover we 
establish a one-to-one correspondence between the $k$ copies of the 
unit circle and the faces $F'$ of $M$ containing $F$, of dimension one 
higher than that of $F$. We then identify the points of the disjoint 
union $\cup F \times (S^1)^k$ as follows.  If $F \subset F'$ and $F'$ 
corresponds to the variable $\theta_i\in S^1$ we identify 
$(x,\theta_1,\ldots,\theta_{i-1},1,\theta_{i+1},\ldots,\theta_{k})\in 
F \times (S^1)^k$ to the point 
$(x,\theta_1,\ldots,\theta_{i-1},\theta_{i+1},\ldots,\theta_{k})\in F' 
\times (S^1)^{k-1}$ (same $x$). The resulting quotient space is by 
definition $\LX$. By construction there exists a continuous map 
$p_{\maL}:\LX \to M$. 
Let $J = S^1 \cup [1, 1+\epsilon)$, for some $\epsilon >0$, with $S^1$
identified with a subset of the complex plane.  Then the space
$\maL(M)$ is locally modelled by $J_\epsilon^k \times \RR^{n-k}$,
above each point of $M$ belonging to an open face of codimension
$k$. Using the space $\maL(M)$, one can describe also the homology of
other complexes associated to $M$.
 
Let $X \subset M$ be a closed subset consisting of a closed union of 
faces of $M$. Denote by $\Omega_c^{k}(M,X)$ the space of compactly 
supported smooth forms on $M$ that vanish to infinite order on $X$. We 
then obtain complexes 
\begin{equation}\label{eq.mou2,1} 
	\cdots\longrightarrow \Omega_c^k(M,X)
	\stackrel{d}{\longrightarrow} \Omega_c^{k+1}(M,X)
	\longrightarrow \cdots
\end{equation} 
and 
\begin{equation}\label{eq.mou2,2} 
	\cdots\longrightarrow \Omega_{{c,\maL}}^k(M,X)
	\stackrel{d}{\longrightarrow} \Omega_{{c,\maL}}^{k+1}(M,X)
	\longrightarrow \cdots,
\end{equation} 
where $\Omega_{{c,\maL}}^k(M,X) := \mOt \Omega_c^k(M,X)$, using the 
same convention as above. We denote by $\cohom_{{\maL}}^*(M,X)$ the 
cohomology of this complex.  When $X = \emptyset$, this recovers the 
old definitions.

\begin{proposition}\label{prop.H*.L}\ 
The homology of the complex \eqref{eq.mou2,1} is the relative (de 
Rham) cohomology group $\cohom_c^*(M,X)$ and the cohomology of the 
complex \eqref{eq.mou2,2} is 
\begin{equation*} 
	\cohom_{\maL}^*(M,X) \cong \cohom_c^*(\LX,p_{\maL}^{-1}(X)) =
	\cohom_{c}^* (\LX \smallsetminus p_{\maL}^{-1}(X)).
\end{equation*} 
Moreover, the morphism on cohomology induced by the inclusion of the 
first complex into the second complex identifies, up to isomorphism, 
with the morphism $p_{\maL}^*$. 
\end{proposition}

\begin{proof}\ 
That the cohomology of the first complex is isomorphic to 
$\cohom_c^*(M,X) :=\cohom_c^*(M \smallsetminus X)$ is of course well 
known.  The second isomorphism follows from the Serre-Leray spectral 
sequence applied to the map $p_{\maL}:\LX \to M$ and the above proposition. 
\end{proof}

In particular, the above theorem gives formula $\cohom_{{ \maL}}^*(M) 
\cong H_c^*(\LX)$ for the Laurent--de Rham cohomology of the manifold 
with corners $M$ (assumed to be $\sigma$-compact).

\section{Homology of complete symbols\label{Sec.gen.HOM}}

Let $\GR$ be a differentiable groupoid with space of units denoted 
$M$, a manifold with corners, in general. We obtain in this section 
some results on the homologies of algebras of complete symbols on 
$\GR$. In particular, we get a complete determination of the periodic 
cyclic cohomology of the algebras of complete symbols $\alge(M,X)$, 
$\alge_0(M,X)$, and $\alge_{\maL}(M,X)$ defined below. 
 
Let $X$ be a locally finite union of closed faces of the manifold with 
corners $M$.  Also, let $\ideal_X \subset \CI(M)$ be the subset of 
functions that vanish to infinite order on $X$. Then we define $\alge 
(M,X) := \ideal_X(\tPS{\infty}/\tPS{-\infty})$ the algebra of complete 
symbols on $\GR$, supported above a compact set of $M$ and vanishing 
to infinite order at $X$.  Similarly, we let 
\begin{equation*} 
\begin{gathered} 
\alge_0 (M,X) := \ideal_X(\tPS{0}/\tPS{-\infty}), \\ 
\alge_{\maL} (M,X) := \mO \ideal_X(\tPS{\infty}/\tPS{-\infty}). 
\end{gathered} 
\end{equation*} 
The algebra $\alge_{\maL}(M,X)$ will be sometimes referred to as the 
algebra of ``Laurent complete symbols vanishing rapidly at $X$.'' 
Note that we have $\alge_{\maL}(M,X) = \alge (M) \mathcal O(M) 
\ideal_X$, too.  If $X = \emptyset$, we do not include it in the 
notation. 
 
Let us recall some notation from \cite{LauterNistor}, in this volume. We 
denote 
by $T_{vert}(\GR)=\cup_{x\in M} T(\GR_x)$ the vertical tangent bundle 
of $\GR$ to the fibers of the domain map $d: \GR \to M$. Then $A(\GR)$ 
can be identified with the restriction of $T_{vert}(\GR)$ to $M$, that 
is, $A(\GR)=\cup_{x\in M} T_x(\GR_x)$. In particular, the smooth 
sections of $A^*(\GR)$ are canonically identified with the right 
$\GR$-invariant smooth sections of $T_{vert}^*(\GR)$.

\begin{proposition}\label{prop.alge}\ 
Suppose $M$, the space of units of $\GR$, is $\sigma$-compact. Then the 
algebras $\alge(M,X)$, $\alge_0(M,X)$, and $\alge_{\maL}(M,X)$ are 
topologically filtered algebras. 
\end{proposition}

\begin{proof}\ 
We shall prove this lemma for $\alge_{\maL}(M,X)$, the other case 
being completely similar, and actually even simpler.  The subspaces 
$F_{n,l}^m \alge$ are defined in the following way. Fix an increasing, 
countable exhaustion $M = \cup K_m$ by compact submanifolds.  For any 
vector bundle $E \to M$, we denote by ${\mathcal S}^n_{rc}(E)$ the 
space of classical symbols supported above some compact set in the 
base $M$. We then choose a quantization map $q : {\mathcal 
S}^n_{rc}(A^*(\GR)) \to \tPS{n}$ as in \cite{NistorWeinsteinXu}, whose 
main property is that it induces a bijection 
\begin{equation*} 
	\mathcal S^n_{rc}(A^*(\GR))/\mathcal S^{n'}_{rc}(A^*(\GR))
	\cong \tPS{n}/\tPS{n'},
\end{equation*} 
for all $n' < n$ (including $n = \infty$ or $n' = -\infty$). Denote by 
$\pi : A^*(\GR) \to M$ the natural projection and by ${\mathcal 
S}^n_{K_m}(\alge^*(\GR))$ the set of symbols with support in 
$\pi^{-1}(K_m)$, then an increasing triple filtration for $\alge$ is 
defined, for $l = 0$ first, by 
\begin{equation*} 
	F_{n,0}^m\alge = \big ( q({\mathcal S}^n_{K_m}(A^*(\GR))) +
	\Psi_{K_m}^{-\infty}(\GR) \big )/ \Psi_{K_m}^{-\infty}(\GR).
\end{equation*} 
 
Choose now a radial completion of $A^*(\GR)$, which is then a 
diffeomorphism of $A^*(\GR)$ onto the interior of the ball bundle 
$B^*(\GR) := \{ \xi \in A^*(\GR), \|\xi \| \le 1 \}$. This identifies 
$S^0_K(A^*(\GR))$ with the subset $\CI_K(B^*(\GR)) \subset 
\CI(B^*(\GR))$ of those smooth functions on $B^*(\GR)$ with support 
above the compact set $K \subset M$. We use this identification to 
define the topology on ${\mathcal S}_{K_m}^0(A^*(\GR))$, which in turn 
gives $F_{0,0}^m\alge$ the induced topology. 
 
To define $F_{n,l}^m\alge$ in general, let ${\mathcal 
S}^{n,l}_{K_m}(A^*(\GR))$ be the space of Laurent symbols of order 
$\le n$, with support in $\pi^{-1}(K_m)$, and only rational 
singularities with total order $\le l$ in each defining function of a 
hyperface of $M$ 
outside $X$ 
(we do not count negative 
orders). Then we define $F_{n,l}^m\alge$ as the image of $q({\mathcal 
S}^{n,l}_{K_m}(A^*(\GR))$, with the induced topology. 
 
This definition is such that $F_{n,l}^m\alge = F_{n,0}^m\alge$, if $l 
< 0$.  The topology on this space is defined similarly.  In this way, 
$F_{n,l}^m\alge$ becomes a closed subset of $F_{n',l'}^{m'}\alge$ 
whenever $n \le n'$, $l \le l'$, and $m \le m'$. Then we endow $\alge$ 
with the strict inductive limit topology. Conditions 1--7 of a 
topologically filtered algebra are then satisfied. 
\end{proof} 
 
\vspace*{1mm} {\em Remark.}\ When $M$ is compact, there is no need to 
consider the additional filtration with respect to $m$. Also, for the 
algebras $\alge_0(M,X)$, there is no need to consider the variable 
$l$, as the filtration is independent of $l$ 
\vspace*{1mm} 
 
We now apply some of the results of the previous sections to the 
algebras $\alge(M,X)$ and $\alge_{\maL}(M,X)$.

In this section, we first concentrate on the cyclic homology, because 
this will lead to a complete determination of the periodic cyclic 
homology of these algebras. We concentrate first on the algebra 
$\alge(M) = \tPS{\infty}/\tPS{-\infty}$ of complete symbols on $\GR$, 
{\em which we shall denote by $\alge$ throughout the rest of this 
section}, for simplicity. The groupoid $\GR$ is arbitrary but fixed in 
what follows. 
 
Fix a metric on $A(\GR)$ and let $P$ be a pseudodifferential operator 
of order one such that $\sigma_1(P) \equiv r$ (modulo lower order 
symbols), where $r \in \CI(A^*(\GR))$ is the distance function to the 
origin. We know that the graded algebra $Gr(\alge)$ associated to 
$\alge$ is commutative. Denote by $S^*(\GR) = S^*(A(\GR))$ the set of 
vectors of length one in $A^*(\GR)$, the dual of the Lie algebroid of 
$\GR$.  Then $Gr(\alge) \simeq \CIc(S^*(\GR)) \otimes \CC[r,r^{-1}]$, 
with grading given by the powers of $r$. 
 
As noted in Section \ref{Sec.HandC}, the tensor products appearing in 
the Hochschild complex are such that $F_k\mathcal 
H_n(\alge)/F_{k-1}\mathcal H_n(\alge)$ is a direct sum of spaces, each 
of which is isomorphic to $\CIc(S^*(\GR) \times S^*(\GR) \times \cdots 
\times S^*(\GR))$, such that the natural map 
\begin{equation*} 
F_k\mathcal H_n(\alge) \to \lim_{\leftarrow}\, F_k\mathcal 
H_n(\alge)/F_{k-j} \mathcal H_n(\alge),\quad j \to \infty, 
\end{equation*} 
is an isomorphism. The same comments are valid for the cyclic complex.

The Hochschild, cyclic, and periodic cyclic homologies of $Gr(\alge)$ 
are identified using a combination of the Hochschild-Kostand-Rosenberg 
(HKR) isomorphism and a result of Connes, which is the analog of the 
HKR-isomorphism for algebras of $C^\infty$-functions.  We denote by 
$\Omega_{rc}^{l}(A^*(\GR) \smallsetminus 0)_d$ the set of 
$l-$differential forms on the manifold $A^*(\GR) \smallsetminus 0$ 
that are positively $d$-homogeneous in the radial direction and whose 
support projects onto a compact subset of $M$. (Here $A^*(\GR) 
\smallsetminus 0$ stands for the dual of $A(\GR)$ with the zero section 
removed.) Then, using the grading of $Gr(\alge)$, we have 
\begin{multline*} 
\Hd_{l}(Gr(\alge))_{d} \cong \Omega_{rc}^{l}(A^*(\GR) 
\smallsetminus 0)_d \cong \Omega_{c}^{l}(S^*(\GR))r^d \oplus 
\Omega_{c}^{l-1}(S^*(\GR))r^{d-1}dr, 
\end{multline*} 
the isomorphism being obtained via the 
Hochschild-Kostant-Rosenberg-Connes map 
\begin{equation*} 
\chi(a_0, \ldots , a_l)=(1/l!)\, a_0 da_1\wedge da_2 \wedge 
\ldots \wedge da_l. 
\end{equation*} 
(We shall often omit the ``wedge'' $\wedge$ in what follows.) 
 
We can identify $\Omega_{c}^{l}(S^*(\GR))r^d \oplus 
\Omega^{l-1}(S^*(\GR))r^{d-1}dr$ with the subspace $\Omega_{rc}^{l} 
(A^*(\GR) \smallsetminus 0)_d \subset \Omega_{rc}^{l}(A^*(\GR) 
\smallsetminus 0)$ consisting of $d$-homogeneous forms, as 
above. Also, we shall sometimes complexify $r$ and restrict it to the 
unit circle, so that $r$ becomes $e^{i\theta}$ and $dr$ becomes $i 
e^{i\theta} d\theta$. Let $\Omega^*_{fin}(S^1)$ be the set of 
polynomial forms on $S^1$, that is, of finite linear combinations of 
$r^l$ and $r^{m}dr$. The complex $\Omega^*_{fin}(S^1)$ has the same 
cohomology as the de Rham complex of all forms on $S^1$, which also 
gives the cohomology of $S^1$ (with complex coefficients). 
 
Then, we see that the complex $(\oplus_m \Omega_{rc}^*(A^*(\GR) 
\smallsetminus 0)_m, d)$, can be identified with the complex 
$\Omega_{c}^*(S^*(\GR)) \otimes \Omega^*_{fin}(S^1)$, and hence it has 
the same cohomology as $S^*(\GR) \times S^1$. Next we use that 
\begin{equation*} 
\begin{gathered} 
\Omega_{rc}^{k+1}(A^*(\GR) \smallsetminus 0)_m/ 
d\Omega_{rc}^k(A^*(\GR) \smallsetminus 0)_m \cong 
\Omega_{c}^{k+1}(S^*(\GR)), \quad \text{ for }\, m \not = 0, 
\;\text{ and }\\ \quad \Omega_{rc}^{k+1}(A^*(\GR) 
\smallsetminus 0)_0/ d\Omega_{rc}^k(A^*(\GR) \smallsetminus 
0)_0 \cong \\ \Omega_{c}^{k+1}(S^*(\GR))/d 
\Omega_{c}^{k}(S^*(\GR)) \oplus 
\Omega_{c}^{k}(S^*(\GR))/d\Omega_{c}^{k-1}(S^*(\GR)). 
\end{gathered} 
\end{equation*} 
The first isomorphism above is given by 
\begin{equation*} 
(\alpha,\beta)\in \Omega_{c}^{k+1}(S^*\GR)r^m \oplus 
\Omega_{c}^{k}(S^*\GR)r^{m-1}dr \mapsto m^{-1}\alpha - d\beta. 
\end{equation*} 
 
\begin{lemma}\label{lemma.EC1}\ 
The $E^1$-term of the spectral sequence associated by Lemma 
\ref{lemma.sp.sq} to the topologically filtered algebra $\alge$ is 
given by $\EC^1_{k,h} \cong \Omega_{c}^{k+h}(S^*(\GR))$ if $k \not = 
0$ and 
\begin{multline*} 
\EC^1_{0,h} \cong \Omega_{c}^{ h}(S^*(\GR))/d 
\Omega_{c}^{h-1}(S^*(\GR)) \oplus \Omega_{c}^{h-1}(S^*(\GR))/d 
\Omega_{c}^{h-2}(S^*(\GR)) \oplus \\ \oplus_{j>0} 
\cohom_{c}^{h-2j}(S^*(\GR) \times S^1). 
\end{multline*} 
Moreover, the periodicity morphism $S : \EC^1_{k,h} \to \EC^1_{k,h-2}$ 
vanishes if $k \not = 0$ and is the natural projection if $k = 0$. 
\end{lemma}

\begin{proof}\ 
We know from Lemma \ref{lemma.sp.sq} that $\EC_{k,h}^1 \simeq 
\Hc_{k+h}(Gr(\alge))_k$. Using again the HKR-isomorphism, we obtain 
that $\Hc_{m}(Gr(\alge))_d$ is isomorphic to 
\begin{equation*} 
\Omega_{rc}^m(A^*(\GR)\smallsetminus 
0)_d/d\Omega_{rc}^{m-1}(A^*(\GR) \smallsetminus 0)_d \oplus 
\oplus_{j > 0} \cohom^{m-2j}(A^*(\GR) \smallsetminus 0)_d, 
\end{equation*} 
and the operator $S: \EC^1_{k,h} \to \EC^1_{k,h-2}$ is the projection 
which sends $\Omega_{rc}^{k+h}(A^*(\GR)\smallsetminus 0)_k/ 
d\Omega_{rc}^{k+h-1}(A^*(\GR) \smallsetminus 0)_k$ to 0, is the 
inclusion 
\begin{equation*} 
\cohom^{k+h-2}(A^*(\GR) \smallsetminus 0)_k \hookrightarrow 
\Omega_{rc}^{k+h-2}(A^*(\GR)\smallsetminus 
0)_k/d\Omega_{rc}^{k+h-3} (A^*(\GR) \smallsetminus 0)_k 
\end{equation*} 
and the identity on the other factors.  The above computation of 
$\Omega_{rc}^m(A^*(\GR)\smallsetminus 0)/ d\Omega_{rc}^{m-1}(A^*(\GR) 
\smallsetminus 0)$ then gives the result. 
\end{proof}

{}From the above lemma we obtain the following immediate consequence.

\begin{corollary}\label{cor.conditions}\ 
We have $\EC^1_{k,h} = 0$ if $k < 0$ and $k + h > \dim S^*(\GR)$. 
\end{corollary}

The following proposition is essential in determining the periodic 
cyclic homology of $\alge$.

\begin{proposition}\label{HC}\ 
If 
$q > \dim S^*(\GR)$, then $\HC_q(\alge) \cong \oplus_{k \in \ZZ} 
\cohom_c^{q - 2k}(S^*(\GR) \times S^1)$ with $S : \HC_{q + 2}(\alge) 
\to \HC_q(\alge)$ also an isomorphism. 
\end{proposition}

\begin{proof}\ 
We shall use Theorem \ref{theorem.conv1} and Lemma 
\ref{lemma.EC1}. For $q > \dim(S^*\GR)+1$, Lemma \ref{lemma.EC1} gives 
\begin{equation*} 
\EC^1_{k,q-k}=0 \text{ for all } k \not = 0,\, \text{ and } 
\EC^1_{0,q} = \oplus _{j \in \ZZ} \cohom_{c}^{q-2j} (S^*\GR 
\times S^1). 
\end{equation*} 
Moreover, the assumptions of Theorem \ref{theorem.conv1} are 
satisfied, by Corollary \ref{cor.conditions}, and this gives 
\begin{equation*} 
\HC_q(\alge) \simeq \oplus_k \EC^1_{k,q-k} \cong \oplus _{j>0} 
\cohom_{c}^{q-2j} (S^*(\GR) \times S^1). 
\end{equation*} 
This gives the result. 
\end{proof}

{}From this we obtain

\begin{theorem}\label{theorem.HP}\ 
The periodic cyclic homology groups of $\alge = 
\tPS{\infty}/\tPS{-\infty}$, the algebra of complete symbols on $\GR$, 
are given by $\Hp_q(\alge) \cong \cohom_c^{[q]}(S^*(\GR)\times S^1)$. 
\end{theorem}

\begin{proof}\ 
Whenever the periodicity operator $S$ of the exact sequence 
\eqref{eq.SBI} is surjective, we have $\Hp_q(\alge) \cong 
\displaystyle{\lim_{\leftarrow}}\, \HC_{q+2j},$ the projective limit 
being taken with respect to the operator $S$.  The conclusion then 
follows from Proposition \ref{HC}. 
\end{proof}

Similarly, we have the following determination of the periodic cyclic 
homology groups of the algebra $\alge_0 = \tPS{0}/\tPS{-\infty}.$

\begin{theorem}\ 
$\Hp_q(\alge_0) \cong \cohom_c^{[q]}(S^*(\GR))$. 
\end{theorem}

\begin{proof}\ 
The proof is essentially the same as for the corresponding result for 
the algebra $\alge$, so we will be brief.  The filtration on $\alge_0$ 
is induced from the filtration on $\alge$.  This and the specific form 
of the $\EC^1$-terms then give 
\begin{equation*} 
\begin{gathered} 
\EC^1_{k,h}(\alge_0)\simeq \EC^1_{k,h}(\alge)\, \text{ if } k 
< 0 ,\, \text{ and } \EC^1_{k,h}(\alge_0)\simeq \{0\}\, \text{ 
if } k > 0,\\ \text{ and } \EC^1_{0,h}(\alge_0) \cong 
\oplus_{j \in \ZZ} \cohom_c^{h - 2j}(S^*(\GR)) \text{ if } h 
\ge \dim S^*(\GR). 
\end{gathered} 
\end{equation*} 
Moreover, the differential $d^1$ is the same as that for $\alge$ if 
$k\leq 0$, but is trivial for $k > 0$. Thus, if $q > \dim(S^*\GR)$, we 
get as above that 
\begin{equation*} 
\forall k\in \Z\smallsetminus \{0\}, \; \EC^1_{k,q-k}(\alge_0) 
\simeq 0. 
\end{equation*} 
Consequently, 
\begin{equation*} 
  \EC^1_{0,q}(\alge_0) \simeq \oplus_{j \in \ZZ} \cohom_c^{q - 
2j}(S^*(\GR)) =: \cohom^{[q]}(S^*(\GR)). 
\end{equation*} 
\end{proof}

A similar approach can be used to treat variants of the algebra 
$\alge$ when the Schwarz symbols of our operators vanish to infinite 
order at certain hyperfaces of $M$. 
 
Let $X$ be a closed union of faces of the manifold with corners $M$. 
Let $\ideal_X \subset \CI(M)$ be the ideal of functions that vanish to 
infinite order on $X$. We now study the algebras $\alge_{\maL}(M,X)$ 
of ``Laurent complete symbols vanishing rapidly at $X$,'' supported 
above a compact subset of $M$, which, we recall, are given by 
\begin{equation*} 
\alge_{\maL}(M,X) := \mathcal O(M) \ideal_X \alge (\GR) = 
\mathcal O(M) \ideal_X (\tPS{\infty}/\tPS{-\infty}), 
\end{equation*} 
and $\alge_{\maL}(\GR) = \alge_{\maL}(\GR,\emptyset)$. 
 
We shall usually denote by $\pi_Y : Y \to M$ the projection associated 
with a typical fibration $Y$ over $M$.  As before, we denote by 
$\Omega_{rc,\maL}^t(A^*(\GR)\smallsetminus 0, 
\pi_{A^*(\GR)\smallsetminus 0}^{-1}(X))$ the space of $t$-differential 
forms on $A^*(\GR)\smallsetminus 0$ that vanish to infinite order in 
the base variable on $\pi_{S^*(\GR)\setminus 0}^{-1}(X)$ and are supported 
above a compact 
set in $M$.  The above results on the cyclic homology 
of the algebras $\alge$ and $\alge_0$ extend then almost right away to 
the algebras $\alge_{\maL}(M,X)$. Recall that $p_{\maL}:\maL(Y) \to Y$ 
is the projection defined in section \ref{Sec.PREM} for a manifold 
with corners $Y$. We summarize results in the following two 
propositions. To simplify notation, we shall denote also by $p_{\maL} 
: \maL(S^*(\GR)) \times S^1 \to M$ the induced projection.

\begin{proposition}\label{hom.AL}\ 
For $q > \dim(S^*(\GR))$, we have 
\begin{equation*} 
	\HC_q(\alge_{\maL}) \cong \oplus_{k \ge 0}\cohom_{c}^{q-2k}
	(\maL(S^*(\GR))\times S^1 \smallsetminus p_{\maL}^{-1}(X) ).
\end{equation*} 
Thus, $\Hp_q(\alge_{\maL})\cong \cohom_{c}^{[q]}(\maL(S^*(\GR))\times 
S^1 \smallsetminus p_{\maL}^{-1}(X)).$ 
\end{proposition}

\begin{proof}\ 
We first prove the analogue of Lemma \ref{lemma.EC1} in our new 
settings: namely, the $E^1$-term of the cyclic spectral sequence 
associated to the topologically filtered algebra $\alge_{\maL}(M,X)$ 
by Lemma \ref{lemma.sp.sq} is given by 
\begin{equation*} 
\EC^1_{k,h} \cong \Omega_{c,\maL}^{k+h} 
(S^*(\GR),\pi_{S^*(\GR)}^{-1}(X)), 
\end{equation*} 
if $k \not = 0$, and otherwise by 
\begin{multline*} 
\EC^1_{0,h} \cong \Omega_{c,\maL}^{ h}(S^*(\GR),\pi_ 
{S^*(\GR)}^{-1}(X))/d 
\Omega_{c,\maL}^{h-1}(S^*(\GR),\pi_{S^*(\GR)}^{-1}(X))\\ 
\oplus 
\Omega_{c,\maL}^{h-1}(S^*(\GR),\pi_{S^*(\GR)}^{-1}(X))/d 
\Omega_ {c,\maL}^{h-2}(S^*(\GR),\pi_{S^*(\GR)}^{-1}(X)) \oplus 
\\ \oplus_{j>0} \cohom_{c}^{h-2j}(\maL(S^*(\GR)) \times 
S^1,p_{\maL}^{-1}(X) ). 
\end{multline*} 
Moreover, the periodicity morphism $S : \EC^1_{k,h} \to \EC^1_{k,h-2}$ 
vanishes if $k \not = 0$ and is the natural projection if $k = 0$. 
 
Indeed, Lemma \ref{lemma.sp.sq}(i) together with a relative version of 
the HKR isomorphism give that 
\begin{multline*} 
\EC^1_{k,h}\cong 
{\frac{\Omega_{rc,\maL}^{k+h}(A^*(\GR)\smallsetminus 0, 
\pi_{A^*(\GR)\smallsetminus 0}^{-1}(X))_k} {d 
\Omega_{rc,\maL}^{k+h-1}(A^*(\GR)\smallsetminus 
0,\pi_{A^*(\GR)\smallsetminus 0} ^{-1}(X))_k}} \oplus \\ 
\oplus_{j>0} \cohom_{rc}^{k+h-2j}(A^*(\GR)\smallsetminus 
0,\pi_{A^*(\GR)\smallsetminus 0} ^{-1}(X))_k. 
\end{multline*} 
Moreover, it follows that the operator $S$ is the projection just 
above.  Now we can compute each term and use the homotopy invariance 
of the relative de Rham cohomology to see that the $k$-homogeneous 
relative cohomology groups $H_{\maL}^{h-2j}(A^*(\GR)\smallsetminus 
0,\pi_{A^*(\GR)\smallsetminus 0}^{-1}(X))_k$ vanish for $k\not =0$ and 
that for $k=0$ we have 
\begin{equation*} 
H_{rc,\maL}^{h-2j}(A^*(\GR)\smallsetminus 
0,\pi_{A^*(\GR)\smallsetminus 0}^{-1}(X))_0 \cong 
H_{c,\maL}^{h-2j}(S^*(\GR)\times S^1, \pi_{S^*(\GR)\times 
S^1}^{-1}(X)). 
\end{equation*} 
On the other hand, we have 
\begin{equation*} 
{\frac{\Omega_{rc,\maL}^{k+h}(A^*(\GR)\smallsetminus 0,\pi_ 
{A^*(\GR)\smallsetminus 0}^{-1}(X))_k} {d 
\Omega_{rc,\maL}^{k+h-1}(A^*(\GR)\smallsetminus 0, 
\pi_{A^*(\GR)\smallsetminus 0}^{-1}(X))_k}} \cong 
\Omega_{c,\maL}^{k+h}(S^*(\GR),\pi_{S^*(\GR)}^{-1}(X)), 
\end{equation*} 
if $k \not = 0$, and, for $k=0$, we have 
\begin{multline*} 
{\frac{\Omega_{rc,\maL}^{k+h}(A^*(\GR)\smallsetminus 0, 
\pi_{A^*(\GR)\smallsetminus 0}^{-1}(X))_0} {d 
\Omega_{rc,\maL}^{k+h-1}(A^*(\GR)\smallsetminus 0, 
\pi_{A^*(\GR)\smallsetminus 0}^{-1}(X))_0}} \\ \cong 
{\frac{\Omega_{c,\maL}^{h}(S^*(\GR),\pi_{S^*(\GR)}^{-1}(X))} 
{d \Omega_{c,\maL}^{h-1}(S^*(\GR),\pi_{S^*(\GR)}^{-1}(X))}} 
\oplus 
{\frac{\Omega_{c,\maL}^{h-1}(S^*(\GR),\pi_{S^*(\GR)}^{-1}(X))} 
{d \Omega_{c,\maL}^{h-2}(S^*(\GR),\pi_{S^*(\GR)}^{-1}(X))}}. 
\end{multline*} 
 
The first assertion is a direct consequence of the above discussion 
and Theorem \ref{theorem.conv2}.  The computation of $\Hp_q$ 
then follows as in the proof of Theorem \ref{theorem.HP}. 
\end{proof}

Proposition \ref{hom.AL} is formulated in such a way that it remains 
true if we drop the subscript $\maL$ (and if we remove $p_{\maL}$). If 
we proceed then as indicated, we obtain the periodic cyclic homology 
of $\alge(M,X)$. 
 
We now take a closer look at the Hochschild homology of $\alge 
=\tPS{\infty}/\tPS{-\infty}$ and of the other related algebras. 
 
We shall use below the Poisson structure of $A^*(\GR)$, which we now 
recall for the benefit of the reader. The natural regular Poisson 
structure of $T_{vert}^*(\GR)$ induces a Poisson structure on 
$A^*(\GR)$. (This is recalled in \cite{NistorWeinsteinXu}, Lemma 7, 
for example).  If $r:\GR \to M$ is the range map, then the image of 
the differential of $r$ restricted to $A(\GR)$ determines a possibly 
singular foliation $S$ on $M$.  On the other hand, the kernel of $r_*$ 
is a family of Lie algebras whose fiber at $x\in M$ is the Lie algebra 
of the Lie group $\GR_x^x$. When the groups $\GR_x^x$ are 
$0$-dimensional, this foliation has no singularities, and $A(\GR)$ 
becomes the tangent bundle to the leaves of $S$. In this case, the 
Poisson structure of $A^*(\GR)\cong T^*S$ is induced from the 
symplectic structures of the leaves. In general, the Poisson structure 
on $A^*(\GR)$ is defined by a two tensor 
\begin{equation*} 
G \in \CI(A^*(\GR), \Lambda^2(T(A^*(\GR)))) 
\end{equation*} 
so that $\{f,g\}=i_G(df\wedge dg).$ 
 
Let $i_G$ be the contraction by $G$. Then we obtain as in 
\cite{Brylinski} a differential 
\begin{equation} 
\delta := i_G\circ d - d \circ i_G :\Omega^k(A^*(\GR)) \to 
\Omega^{k-1}(A^*(\GR)), 
\end{equation} 
 
Explicitly, for any open subset $V \subset A^*(\GR)$ of the form 
$V\cong [0,1)^l\times \R^{n-l}$ and any $(f_0,...,f_k)\in 
\CI(V)^{k+1}$, the differential $\delta$ is given locally by the 
formula 
\begin{multline*} 
\delta(f_0 d f_1\wedge d f_2 \wedge \ldots \wedge d f_k) = 
   \sum_{1\leq j\leq k} (-1)^{j+1} \{f_0,f_j\}d f_1\wedge \ldots 
   \wedge {\widehat{d f_j}}\wedge \ldots \wedge d f_k \\ + 
   \sum_{1\leq i<j\leq k} (-1)^{i+j}f_0d \{f_i,f_j\} \wedge d f_1 
   \wedge \ldots \wedge {\widehat{d f_i}} \wedge \ldots \wedge 
   {\widehat{d f_j}} \wedge \ldots \wedge d f_k. 
\end{multline*} 
This formula is valid also when $M$ has corners.  It is easy to check 
that the differential $\delta$ is homogeneous of degree $-1$ with 
respect to the action of $\R^*_+$ on $A^*(\GR)\smallsetminus 0$ and 
hence maps 
\begin{equation*} 
\delta : \Omega^p(A^*(\GR)\smallsetminus 0)_{d} \to 
\Omega^{p-1} (A^*(\GR)\smallsetminus 0)_{d-1}. 
\end{equation*} 
Moreover, $\delta$ preserves the support, so it maps the space 
$\Omega_{rc}^p(A^*(\GR)\smallsetminus 0)_{d}$ of $d$-homogeneous forms 
supported above a compact set of $M$ to 
$\Omega_{rc}^{p-1}(A^*(\GR)\smallsetminus 0)_{d-1}.$

\begin{definition}\label{poisson}\ 
The $d$-homogeneous Poisson $p$-homology space $\cohom_p^{\delta} 
(A^*(\GR)\smallsetminus 0)_d$ of the conic Poisson manifold with 
corners $A^*(\GR)\smallsetminus 0$ is defined by 
\begin{equation*} 
\cohom_p^{\delta}(A^*(\GR)\smallsetminus 0)_d := {\frac{\ker 
(\delta : \Omega_{rc}^p(A^*(\GR)\smallsetminus 0)_{d} \to 
\Omega_{rc}^{p-1}(A^*(\GR)\smallsetminus 0)_{d-1})} 
{\delta(\Omega_{rc}^{p+1}(A^*(\GR)\smallsetminus 0)_{d+1}}}. 
\end{equation*} 
\end{definition}

Let $\Omega_{rc,\maL}^p(A^*(\GR)\smallsetminus 0)_d$ be the set of 
Laurent differential $p$-forms which are $d$-homogeneous under the 
radial action of $\R^*_+$ and such that their support projects onto a 
compact subset of $M$.  Then it can be checked that $\delta$ is well 
defined on $\Omega_{{rc,\maL}}^*(A^*(\GR)\smallsetminus 0)_d$ and is 
given by the same formula. In fact the covector $G$ is here tangent to 
the faces because the groupoid $\GR$ is differentiable. We thus have 
the following Laurent complexes indexed by $p \in \Z$ 
\begin{equation*} 
0\to {\mathcal P}^{p,n - p} \stackrel{\delta}{\rightarrow} 
{\mathcal P}^{p,n - p -1}\stackrel{\delta} {\rightarrow} 
...\stackrel{\delta}{\rightarrow} {\mathcal P}^{p,-p}\to 0, 
\end{equation*} 
where ${\mathcal 
P}^{p,d}=\Omega_{{rc,\maL}}^{p+d}(A^*(\GR)\smallsetminus 0)_d$.  The 
homology of the resulting complex is, by definition, the 
Laurent-Poisson homogeneous homology of $A^*(\GR)\smallsetminus 
0$. This homology will be denoted by 
\begin{equation*} 
\cohom_{\maL,p+d}^{\delta}(A^*(\GR)\smallsetminus 0)_d := 
{\frac{\ker(\delta:{\mathcal P}^{p,d}\to {\mathcal 
P}^{p,d-1})} {\delta(\mathcal P^{p,d+1})}}. 
\end{equation*}

We begin again with a result on the algebra $\alge = 
\tPS{\infty}/\tPS{-\infty}$.

\begin{proposition}\label{prop.E1}\ 
The algebra $\alge$ is $H$-unital.  Let 
\begin{equation*}
	\chi: \HH_{l}(Gr(\alge))_d \rightarrow
	\Omega_{rc}^{l}(A^*(\GR) \smallsetminus 0)_d
\end{equation*}
be the HKR-isomorphism, and let $d_1 : \EH^1_{k,h}\rightarrow
\EH^1_{k-1,h}$ be the first differential of the spectral sequence
associated to $\alge$ by Lemma \ref{lemma.sp.sq}.  Then $\chi \circ
d_1 \circ \chi^{-1}= -\sqrt{-1} \delta$, and hence $\EH_{k,h}^2 \simeq
\cohom^\delta_{k+h}(A^*(\GR)\smallsetminus 0)_k$. 
\end{proposition}

\begin{proof}\ 
The $H$-unitality follows from Propositions \ref{prop.Hunital} and 
\ref{prop.alge}.  For the second part, we proceed essentially as in 
\cite{BrylinskiGetzler}, Theorem 1. Choose an anti-symmetric tensor in 
the last $m$-variables 
\begin{equation*} 
	\eta = \sum \operatorname{sign}(\sigma) f_0 \otimes 
	f_{\sigma(1)} \otimes \ldots \otimes f_{\sigma(m)}, 
\end{equation*} 
with $f_j \in S_c^\infty(A^*(\GR))$. We denote by 
\begin{equation*} 
	q(\eta) = \sum \operatorname{sign}(\sigma) q(f_0) \otimes 
	q(f_{\sigma(1)}) \otimes \ldots \otimes q(f_{\sigma(m)}) 
\end{equation*} 
the {\em quantization} of $\eta$. Let $k = \deg f_0 + \ldots + \deg 
f_m$ be the total degree.  Because 
\begin{equation*}
	[q(a), q(b)] = -\sqrt{-1} q(\{a,b\}) + \ldots ,
\end{equation*}
where the dots represent terms of order at most $\deg a + \deg b - 2$,
the quantity $b \circ q(\eta)$ is of total order at most $k-1$ and
hence, modulo terms of order $k-2$, $\chi
\circ b \circ q(\eta))$ is easily checked to be exactly 
$\delta(\eta)$. 
\end{proof}

We think that the above spectral sequence is always convergent and 
that it actually degenerates at $E^2$.  We prove this, for instance, 
for differentiable groupoids associated with manifolds with corners, 
see in Section \ref{Sec.HH}.  This spectral sequence is also studied 
in \cite{BenameurNistor} for foliations. In \cite{CrainicMoerdijk} 
it is proved that the groupoid algebras $\tPS{-\infty}$ are 
$H$-unital. 
 
In the same way we can generalize Proposition \ref{prop.E1} and state:

\begin{proposition}\label{HH.AL}\ 
The algebra $\alge_{\maL}(M,X)$ is $H$-unital.  The $\EH^2$ term of 
the spectral sequence associated by Lemma \ref{lemma.sp.sq} to 
$\alge_{\maL}(M,X)$ is given by 
\begin{equation*} 
\EH_{k,h}^2 \simeq \cohom^\delta_{{\maL},k+h}( 
(A^*(\GR)\smallsetminus 0) \smallsetminus 
\pi_{A^*(\GR)\smallsetminus 0}^{-1}(X))_k. 
\end{equation*} 
\end{proposition}

\begin{proof}\ 
The $H$-unitality follows from Propositions \ref{prop.Hunital} and 
\ref{prop.alge}, as above. 
 
To prove the rest of this proposition, we can either use the same 
method as the one used to prove Proposition \ref{prop.E1}, or we can 
argue that this proposition actually follows from Proposition 
\ref{prop.E1}. 
\end{proof}

\section{Hochschild homology for manifolds with corners \label{Sec.HH}}

We now restrict our study to groupoids $\GR$ in a particular class. 
Thus, we shall assume throughout the rest of this paper that the 
anchor map $\varrho: A(\GR) \to TM$ is such that the induced map 
$\varrho_* : \Gamma(A(\GR)) \to \Gamma(TM)$ of $\CI(M)$-modules 
becomes an isomorphism after tensoring with $\mO$. That is, we assume 
that 
\begin{equation} \label{eq.bound-fibr} 
	\varrho_*:\mOt \Gamma(A(\GR)) \cong \mOt \Gamma(TM). 
\end{equation} 
We then say that $A(\GR)$ is {\em rationally isomorphic to $TM$}.  If 
that is the case, it also follows that $\varrho$ induces an 
isomorphism of $A(\GR)$ to $TM$ on the interior $M_0$ of $M$: 
\begin{equation*} 
	A(\GR)\vert_{M_0} \cong TM_0. 
\end{equation*} 
Let $\mO$ be the ring of smooth functions on the interior $M$ with 
only Laurent singularities at the corners, as in Section 
\ref{Sec.PREM}.  We proceed in this section to the 
computation of the Hochschild homology groups of $\alge : = 
\tPS{\infty}/\tPS{-\infty}$, the algebra of complete symbols on $\GR$ 
and of various other related algebras. 

Throughout this and the following sections, $n$ denotes the dimension
of $M$, the space of units of $\GR$.
 
Let $X \subset M$ be a closed union of faces of $M$ and denote by 
$\ideal_X \subset \CI(M)$ the ideal of functions that vanish to 
infinite order on $X$. Recall that we have defined the algebra 
$\alge_{\maL}(M,X)$ of ``Laurent complete symbols vanishing rapidly at 
$X$'' by 
\begin{equation*} 
\alge_{\maL}(M,X) := \mathcal O(M) \ideal_X \alge = \mathcal 
O(M) \ideal_X (\tPS{\infty}/\tPS{-\infty}), 
\end{equation*} 
($\alge_{\maL}(M) = \alge_{\maL}(M,\emptyset)$). 
 
We keep the notations of Section \ref{Sec.gen.HOM} so that $\pi_{S^*M 
\times S^1} : S^*M \times S^1 \to M$ is the natural projection.  Also, 
note that $S^*(\GR) \cong S^*M$, but not canonically. The homotopy 
class of this homeomorphism is canonical though. Also, $p_{\maL}$ 
denotes the structural map $\maL(S^*(\GR)) \times S^1 \to M$.

\begin{theorem}\label{theorem.alg}\ 
Let $X \subset M$ be a union of hyperfaces of $M$ and assume that 
$A(\GR)$ is rationally isomorphic to $TM$ (that is, it satisfies the 
condition \eqref{eq.bound-fibr} above). Then we have 
\begin{equation*} 
\Hd_*(\alge_{\maL}(M,X))\simeq \cohom_{c}^{2n-*}( \maL(S^* M) 
\times S^1 \smallsetminus p_{\maL}^{-1}(X) ) 
\end{equation*} 
\end{theorem}

\begin{proof}\ 
The computation of $\Hd_*(\alge_{\maL}(M,X))$ proceeds exactly as in 
the case of manifolds with boundary. {}From Proposition \ref{HH.AL} we 
know that the $E^2$ term of the spectral sequence associated with the 
order filtration on the Hochschild complex is given by 
\begin{equation*} 
\EH^2_{k,h}\cong \cohom_{\maL,k+h}^\delta 
(A^*(\GR)\smallsetminus 0, \pi_{A^*(\GR)\smallsetminus 
0}^{-1}(X))_k. 
\end{equation*} 
 
In our case, 
\begin{equation*} 
\Omega_{rc,\maL}^{k+h}(A^*(\GR)\smallsetminus 0,\pi_{A^*(\GR) 
\smallsetminus 0}^{-1}(X))_k \cong \Omega_{rc,\maL}^{k+h}( 
T^*M\smallsetminus 0, \pi_{ T^*M \smallsetminus 0}^{-1}(X))_k, 
\end{equation*} 
which is the main reason why this result is true, and explains the 
thinking behind the condition \eqref{eq.bound-fibr}. Thus the Poisson 
structure on $A^*(\GR)\smallsetminus 0$ can be related to the natural 
Poisson structure on the cotangent bundle $T^*M \smallsetminus 0$. In 
particular, the restriction of the given Poisson structure to the 
interior of $A^*(\GR)\smallsetminus 0$ is a symplectic structure. 
 
We now recall the definition of the symplectic Hodge operator $*_G$. 
Let $\omega$ be the two form on the interior of $M$ that defines the 
symplectic form. The form $\omega$ has only Laurent singularities at 
the hyperfaces due to our assumptions on $A^*(\GR)$. This then defines 
the $*_G$-operator by the formula \cite{Brylinski} 
\begin{equation*} 
*_G : \Omega_{rc,\maL}^t(T^*M) \to 
\Omega_{rc,\maL}^{2n-t}(T^*M), \quad \alpha \wedge 
*_G(\alpha')=\wedge^n(G)(\alpha,\alpha') \omega^n/n!. 
\end{equation*} 
Then $*_G^2=1$ and $*_G \circ \delta \circ *_G = (-1)^k d$ on the 
interior of $M$, by \cite{Brylinski}, and hence also on 
$\Omega^k_{rc,\maL}(T^*M \smallsetminus 0)$ because the forms in the 
latter space are determined by their restriction to the interior of 
$M$. 
 
Now, keeping in mind that by construction the operator $*_G$ maps 
$k$-homogeneous $q$-forms into $k+n-q$-homogeneous $(2n-q)$-forms, we 
get 
\begin{multline*} 
\EH^2_{k,h} \cong H^{\delta}_{\maL,k+h}(A^*(\GR)\smallsetminus 
0, \pi_{A^*(\GR)\smallsetminus 0}^{-1}(X))_k \\ \cong 
H^{2n-k-h}_{\maL}(T^*M\smallsetminus 0,\pi_{T^*M 
\smallsetminus 0}^{-1}(X))_{n-h}. 
\end{multline*} 
 
Using Proposition \ref{prop.H*.L} and the relative homotopy invariance 
of the relative De Rham cohomology, we see that the only non-zero 
homogeneous component corresponds to $h = n$, and then that it is 
given by 
\begin{equation*} 
\EH^2_{k,n}\cong \cohom^{n-k}_{\maL}(S^* M \times S^1, 
\pi_{S^*M \times S^1}^{-1}(X)). 
\end{equation*} 
{}From the fact that $\EH^2_{k,h}=0$, if $h\not = n$, we see that 
$d_r=0$, for all $r \ge 2$, and hence we deduce the result.  The 
convergence of the spectral sequence to $\HH_*(\alge_{\maL}(M,X))$ 
follows from Theorem \ref{theorem.conv2}. 
 
By replacing $k$ with $k - n$ we obtain the stated formula for 
$\Hd_k(\alge_{\maL}(M,X))$. 
\end{proof}

Let us look now in more detail at one of the constructions used in the 
above proof in a particular case.  Suppose then that the sections of 
the vector bundle $A(\GR)$ at the hyperface $H = \{ x = 0 \}$ are 
generated by $\pa_y$ and $x^{c_H}\pa_{x}$. Here $c_H$ are integers 
$\ge 1$ associated to each face and the $y$'s are coordinates on $H$, 
using a tubular neighborhood of $H$. (This definition depends on a 
choice of defining functions for the hyperfaces of $M$.)  Let us call 
this calculus the $c_n$-calculus. 
 
Then the 2-covector is a covector on the manifold with corners 
$A^*(\GR) \smallsetminus 0$ which is given in any open subset of $M$ 
diffeomorphic to one of the form $[0,1)^k\times \R^{n-k}$ (with 
defining functions $(x_1, \ldots, x_k)$) by 
\begin{equation*} 
G = \sum_{k+1\leq j\leq n} \partial_{\xi_j} 
\wedge\partial_{y_j} + \sum_{1\leq j\leq k} x_j^{c_j} 
\partial_{\xi_j} \wedge \partial_{x_j} . 
\end{equation*} 
Here the convention of sign is that of \cite{Brylinski}.  The 
symplectic form is then given by 
\begin{equation*} 
\omega=\sum_{k+1\leq j\leq n} dy_j \wedge d\xi_j + \sum_{1\leq 
j\leq k} x_j^{-c_j} dx_j \wedge d\xi_j. 
\end{equation*} 
 
These two formul\ae\ can be checked in the case $[0,1)\times \R$ and 
globalized just like in the case of the usual cusp calculus. We have 
for instance in the case of a manifold with boundary given by $\{x = 
0\}$ that $*_G(f)=\pm (f/x^c)dx\wedge d\xi$ while $*_G(f dx\wedge 
d\xi)=\pm f x^c$.

\section{Applications\label{Sec.Applications}}

We keep the assumptions of the previous section.  The same methods as 
the ones in the previous section can be used also to determine the 
homology of some quotient algebras. Let $X$ be a closed union of faces 
of $M$.  First we prove a result that allows us to reduce the 
computation of the homologies of the relative algebras 
$\alge_{\maL}(M,X)$ to the computation of those of an ``absolute 
algebra.'' 
 
More precisely, denote by $\alge_{\maL}(M\smallsetminus X)$ the 
algebra of Laurent complete symbols on $\GR\vert_{M\smallsetminus X}$, 
the groupoid obtained by restricting $\GR$ to the invariant subset 
$M\smallsetminus X$.  Also, let $\maL(S^*M)$ be the space introduced 
at the end of Section \ref{Sec.PREM} (with $M$ in place of $S^*M$), 
and denote by $p_{\maL}:\maL(S^*M) \times S^1 \to M$ the resulting 
natural map.  Then we have the following excision result.

\begin{proposition}\ 
The inclusion $\alge_{\maL}(M\smallsetminus X) \to \alge_{\maL}(M,X)$ 
induces isomorphisms in Hochschild, cyclic, and periodic cyclic 
homologies. 
\end{proposition}

\begin{proof}\ 
By a standard argument using the 'SBI'-exact sequence, it is enough to 
prove this statement for Hochschild homology. 
 
Consider then the spectral sequences associated to the two algebras 
and the order filtration on their Hochschild complexes by Lemma 
\ref{lemma.sp.sq}. The induced map is then an isomorphism at the 
$E^2$-term, by Theorem \ref{theorem.alg}, because 
\begin{equation*} 
\maL(S^*M) \times S^1 \smallsetminus p_{\maL}^{-1}(X) = 
\maL(S^*(M \smallsetminus X)) \times S^1. 
\end{equation*} 
\end{proof}

\begin{proposition}\label{prop.comp.exact}\ 
The quotient $\alge := \alge_{\maL}(M)/ \alge_{\maL}(M,X)$ is a 
topologically filtered algebra with 
\begin{equation*} 
	\Hd_q(\alge) \cong \cohom_c^{2n-q}(p_{\maL}^{-1}(X)). 
\end{equation*} 
Moreover, the exact sequence in Hochschild homology associated to the 
exact sequence $0 \to \alge_{\maL}(M,X) \to \alge_{\maL}(M) \to \alge 
\to 0$ is naturally isomorphic to the long exact sequence in 
cohomology associated to the pair $(\maL(S^*M) \times S^1, 
p_{\maL}^{-1}(X))$. 
\end{proposition}

\begin{proof}\ 
Only the last statement is not similar to some other proofs in the 
previous sections. Choose an open neighborhood $U$ of $X$ in $M$ such 
that $X$ is a deformation retract of $U$. Also, let $\phi$ be a smooth 
function supported in $U$ which is one in some smaller neighborhood of 
$X$. This data gives rise to a lifting of any smooth function on 
$\pi^{-1}(X)$ to a function on $S^*M$ with support in $U$, and hence 
also to a linear lift $\psi:\mathcal H_*(\alge) \to \mathcal 
H_*(\alge_{\maL}(M))$. The computation of $[b,\psi]$ identifies the 
boundary map in the Hochschild cohomology spectral sequence and gives 
the result. 
\end{proof}

It is interesting to notice that as a consequence of our computations 
we obtain that the traces live only on the minimal faces. Obviously, 
the faces of $M$ that are manifolds without corners are exactly the 
ones that are minimal with respect to inclusion in the set of faces of 
$M$. Let then $Y$ be the union of all faces of $M$ that have no 
corners (and no boundary).

\begin{proposition}\ 
Let $\tau \in \Hd^0(\alge_{\maL}(M))$ be a trace. Then it is induced 
from a trace on $\alge_{\maL}(M)/ \alge_{\maL}(M,Y)$. 
\end{proposition}

\begin{proof}\ Let $n$ be the dimension of $M$, as before. The dual 
of Proposition \ref{prop.comp.exact} holds true for Hochschild
cohomology. This implies that the inclusion
\begin{equation*}
	\Hd_0(\alge_{\maL}(M)/ \alge_{\maL}(M,Y)) \to
	\Hd_0(\alge_{\maL}(M))
\end{equation*}
is isomorphic (in the sense of the above Proposition) to the dual of
the map $\cohom_c^{2n}(\maL(M)) \to \cohom_c^{2n}(p_{\maL}^{-1}(Y))$.
But this last map is checked to be an isomorphism.
%
%
\end{proof}

The above proof gives as a corollary then the ``number'' of residue
traces on $M$.

\begin{corollary}\ 
The dimension of the space of traces of $\alge_{\maL}(M)$ is the 
number of minimal faces of $M$. 
\end{corollary}

The above two results can be proved directly, without using
singular cohomology groups.


We continue to assume that $A(\GR)$ is rationally isomorphic to $TM$
(that is, it satisfies Equation \eqref{eq.bound-fibr}), as in the
previous section. For the algebras that we have considered, the cyclic
homology can be computed directly in terms of the Hochschild homology
as we shall see below. Let $X \subset M$ be a union of closed faces of
$M$.
 
Consider Connes' `SBI'-exact sequence associated to the algebra 
$\alge_{\maL}(M,X)$ (see Equation \eqref{eq.SBI}).

\begin{theorem}\ 
We have $B=0$ in the `SBI'-exact sequence of the algebra 
$\alge_{\maL}(M,X)$, and hence 
\begin{equation} 
	\Hc_m(\alge_{\maL}(M,X)) = \oplus_{k \geq 0}
	\Hd_{m-2k}(\alge_{\maL}(M,X)).
\end{equation} 
\end{theorem}

\begin{proof}  \ 
We will proceed by induction on $m$ to show that the morphism 
\begin{equation*} 
	B: \Hc_{m-1}(\alge_{\maL}(M,X)) \to \Hd_m(\alge_{\maL}(M,X))
\end{equation*} 
is zero for any $m$. 
 
Let $n$ denote the dimension of the manifold $M$, as above.  The 
statement of the theorem is obviously true for $m=0$.  The proof of 
Theorem \ref{theorem.alg}, show that the groups $\Hd_q 
(\alge_{\maL}(M,X))$ are generated by elements of order $-n+q$ with 
respect to the degree filtration, and that all cycles of order less 
that $-n+q$ are boundaries. This is a direct consequence of the 
computation of the $E^1$ terms of the spectral sequences associated to 
the degree filtration. 
 
Assuming now the statement to be true for all values less than $m$, we 
obtain from the `SBI' exact sequence that the groups 
$\Hc_{m-1}(\alge_{\maL}(M,X))$ are isomorphic to 
$\Hc_{m-1}(\alge_{\maL}(M,X)) \simeq \oplus_{k \geq 0} 
\Hd_{m-2k-1}(\alge_{\maL}(M,X))$. This shows that the groups 
$\Hc_{m-1}(\alge_{\maL}(M,X))$ are generated by elements of order at 
most $-n + m -1$ in the degree filtration. It follows that they map to 
elements of order less than $-n+m$ in $\Hd_m(\alge_{\maL}(M,X))$, and 
hence they vanish by the above remark. 
 
The last statement follows directly from the vanishing of $B$ in the 
Connes' exact sequence. 
\end{proof}

%
%

\appendix 
 
\section{Projective limits} 
 
In this appendix we recall some well known facts about projective (or 
inverse) limits and the homology of projective limits, using this also 
as an opportunity to fix notation. An important role is played by 
$\lim^1 = \displaystyle{\lim_{\leftarrow}}^1$, the first (and only) 
derived functor of the projective limit functor $ 
\displaystyle{\lim_{\leftarrow}}$. 
 
Let $\phi_n : V_{n+1} \to V_{n}$, $n \in \NN$, be an inverse system of 
vector spaces. Define $F : \prod V_n \to \prod V_n$ by $F(v_k) = (v_k 
- \phi_k(v_{k+1}))$. Then $\displaystyle{\lim_{\leftarrow}}\, V_n$ is 
the kernel of $F$ and $\displaystyle{\lim_{\leftarrow}}^1$ is the 
cokernel of $F$, by definition. 
 
Suppose now that $V_n$ are complexes of vector spaces and the maps 
$\phi_n$ are surjective. By writing the homology long exact sequence 
associated to the short exact sequence of complexes 
\begin{equation*} 
	0 \to \displaystyle{\lim_{\leftarrow}}\, V_n \to \prod V_n
	\stackrel{F}{\longrightarrow} \prod V_n \to 0,
\end{equation*} 
we obtain the following well known exact sequence:

\begin{lemma}\label{A.exact.limp}\ 
If the maps $\phi_n : V_{n+1} \to V_n$ are surjective, then the 
homology of the inverse limit satisfies 
\begin{equation*} 
	0 \to \displaystyle{\lim_{\leftarrow}}^1 \cohom_{q+1} (V_n)
	\to \cohom_{q} (\displaystyle{\lim_{\leftarrow}}\, V_n) \to
	\displaystyle{\lim_{\leftarrow}}\, \cohom_q (V_n) \to 0.
\end{equation*} 
\end{lemma}

If $B_n \subset A_n$ is a subspace preserved by $\phi_{n}: A_{n+1} \to 
A_n$ and $C_n = A_n/B_n$, then we have an exact sequence 
\begin{equation*} 
0 \to \prod B_n \to \prod A_n \to \prod C_n \to 0 
\end{equation*} 
and $F$ is an endomorphism of this exact sequence. The ker-coker lemma 
(see \cite{AtiyahMacDonald}, for example) then gives the following 
exact sequence:

\begin{lemma}\label{A.exlim1}\ 
If $A_n$, $B_n$, and $C_n$ are as above, then we have the exact 
sequence 
\begin{equation*} 
0 \to \displaystyle{\lim_{\leftarrow}}\, B_n \to 
\displaystyle{\lim_{\leftarrow}}\, A_n \to 
\displaystyle{\lim_{\leftarrow}} C_n \to 
\displaystyle{\lim_{\leftarrow}}^1 B_n \to 
\displaystyle{\lim_{\leftarrow}}^1 A_n \to 
\displaystyle{\lim_{\leftarrow}}^1 C_n \to 0. 
\end{equation*} 
\end{lemma}

%
%

\end{document}